\documentclass[12pt]{amsart}

\usepackage{amssymb}
\usepackage{amscd}

\newcommand{\bt}{\begin{theorem}}
\newcommand{\et}{\end{theorem}}
\newcommand{\bp}{\begin{proposition}}
\newcommand{\ep}{\end{proposition}}
\newcommand{\bq}{\begin{question}}
\newcommand{\eq}{\end{question}}
\newcommand{\bl}{\begin{lemma}}
\newcommand{\el}{\end{lemma}}
\newcommand{\br}{\begin{result}}
\newcommand{\er}{\end{result}}
\newcommand{\be}{\begin{equation}}
\newcommand{\ee}{\end{equation}}
\newcommand{\bc}{\begin{corollary}}
\newcommand{\ec}{\end{corollary}}

\newtheorem{theorem}{Theorem}[section]
\newtheorem{corollary}[theorem]{Corollary}
\newtheorem{lemma}[theorem]{Lemma}
\newtheorem{proposition}[theorem]{Proposition}
\newtheorem{result}[theorem]{Result}

\newcommand{\N}{\mathbb{N}}

\newcommand{\Z}{\mathbb{Z}}

\newcommand{\cB}{\mathcal{B}}

\newcommand{\cD}{\mathcal{D}}
\newcommand{\cH}{\mathcal{H}}

\newcommand{\cJ}{\mathcal{J}}
\newcommand{\cK}{\mathcal{K}}
\newcommand{\cL}{\mathcal{L}}

\newcommand{\cO}{\mathcal{O}}
\newcommand{\cR}{\mathcal{R}}

\newcommand{\cY}{\mathcal{Y}}

\newcommand{\epr}{\hspace{\fill}$\Box$}

\newcommand{\med}{\medskip}
\newcommand{\sm}{\smallskip}

\begin{document}

\noindent {\em Semigroup Forum} (2009) {\bf 78}, 310 -- 325\\
DOI 10.1007/s00233-008-9101-5
\vspace{0.05in}\\
{\em arXiv version}: fonts, pagination and layout differ from the SF published version; moreover, the published version contains two typos which have been corrected here -- in the paragraph following Remark 3, it is now correctly written that $(ab)a^n=a^n$ and $b^n(ab)=b^n$, and in the paragraph following Lemma 2.6 the word `bisimple' has been inserted on line 4  between the words `any' and `orthodox'. 
\vspace{0.1in}\\
\title{Bisimple monogenic orthodox semigroups}
\author[Simon M. Goberstein]{Simon M. Goberstein*}\thanks{*Partially supported by a CSU research grant.
\vspace{0.04in}\\
\indent Department of Mathematics and Statistics,
California State University, Chico, CA 95929-0525, U.S.A.\\
\indent e-mail: sgoberstein@csuchico.edu}
\begin{abstract}
We give a complete description of the structure of all bisimple orthodox semigroups generated by two mutually inverse elements. 
\vspace{0.1in}
\end{abstract}
\maketitle

\font\caps=cmcsc10 scaled \magstep1   

\section*{Introduction}
\med 
It is well known that cyclic groups and the bicyclic semigroup are the only bisimple monogenic inverse semigroups. On the other hand, the class of bisimple orthodox semigroups generated by two mutually inverse elements is substantially more diverse. The purpose of this paper is to determine the structure of all semigroups of that class using only some standard facts about orthodox semigroups. Partial results in this direction, based on the study of the lattice of congruences on the free orthodox semigroup generated by a pair of mutually inverse elements, can be found in \cite{key8} -- see Section 1 below for more details. 

Let $S$ be an orthodox semigroup generated by two mutually inverse elements $a$ and $b$. In Section 2, we determine the structure of all such bisimple semigroups $S$ when $a$ (and therefore $b$) does not belong to any subgroup of $S$ (see Theorem \ref{28}, which is the main result of the paper). In Section 3, we describe the structure of $S$ when $a$ and $b$ are group elements of $S$ (in which case, $S$ is necessarily bisimple). In the sequel to this article, we will apply its results to the study of lattice isomorphisms of bisimple orthodox semigroups generated by a pair of mutually inverse elements. 
\def\bfseries{\normalsize\caps}
\section{Preliminaries}
\med
Let $S$ be a semigroup. We say that $x\in S$ is a {\em group element} of $S$ if it belongs to some subgroup of $S$; otherwise $x$ is a {\em nongroup element} of $S$. The set of nongroup elements of $S$ will be denoted by $N_S$, and the set of idempotents of $S$ by $E_S$. If $x=yz$ for some $x, y, z\in S$, according to standard terminology, $y$ is a {\em left} and $z$ a {\em right divisor} of $x$. As usual, $\langle X\rangle$ stands for the subsemigroup of $S$ generated by its subset $X$, and $\langle x\rangle$ for the cyclic subsemigroup of $S$ generated by $x\in S$. The order of an element $x$ of $S$ will be denoted by $o(x)$; if $x$ has infinite order, we will write $o(x)=\infty$. If $w=w(x_1,\ldots,x_n)$ is a word in the alphabet $\{x_1,\ldots,x_n\}\subseteq S$, we will say shortly that $w$ is a word in $x_1,\ldots, x_n$, and if no confusion is likely, $w$ will be identified with its value in $S$. For any $x\in S$, we define $x^0$ to be the identity of the semigroup $S^1$, so $x^0y=yx^0=y$ for all $y\in S$ (exceptions of this agreement will occur when $S$ contains a subsemigroup $U$ with an identity $e$ and it will be convenient to put $x^0=e$ for $x\in U$; all such situations will be explicitly identified). Recall that $\cR$ is a left and $\cL$ a right congruence on $S$. These and other standard results concerning Green's relations on $S$ will be used without reference. As shown by Hall \cite[Result 9]{key12}, if $U$ is a regular subsemigroup of $S$, then using superscripts to distinguish Green's relations on $U$ from those on $S$, we have $\cK^U=\cK^S\cap(U\times U)$ for $\cK\in\{\cL, \cR,\cH\}$. This result will also be applied without mention. If $\cH$ is the identity relation on $S$, it is common to say that $S$ is {\em combinatorial}, so a regular semigroup is combinatorial if and only if it has no nontrivial subgroups. 

The {\em bicyclic semigroup} $\cB(a,b)$ is usually defined as a semigroup with identity $1$ generated by the two-element set $\{a,b\}$ and given by one defining relation $ab=1$. Disposing of the identity in the definition of $\cB(a,b)$, one can define it as a semigroup given by the following presentation: $\cB(a,b)=\langle a, b\;|\;aba=a, bab=b, a^2b=a, ab^2=b\rangle$. The structure of $\cB(a,b)$ is well known -- it is a combinatorial bisimple inverse semigroup, each of its elements has a unique representation in the form $b^ma^n$ where $m$ and $n$ are nonnegative integers (and $a^0=b^0=ab$), the semilattice of idempotents of $\cB(a,b)$ is a chain: $ab>ba>b^2a^2>\cdots$, and for each $b^ma^n\in\cB(a,b)$, we have $R_{b^ma^n}=\{b^ma^l\,:\,l\geq 0\}$ and $L_{b^ma^n}=\{b^ka^n\,:\,k\geq 0\}$ (see \cite[Lemma 1.31 and Theorem 2.53]{key5}). These facts will be used below without reference.
 
Let $S$ be a semigroup. Following Clifford \cite{key3}, we write $x\perp y$ to indicate that $xyx=x$ and $yxy=y$ for $x, y\in S$, and the phrase ``$x\perp y$ in $S$'' will be used to express briefly the fact that $x,y\in S$ and $x\perp y$. As usual, the set of all inverses of $x\in S$ will be denoted by $V_S(x)$, sometimes without the subscript $S$ if no confusion is likely. Thus $y\in V_S(x)$ precisely when $x\perp y$ in $S$. An {\em orthodox semigroup} is a regular semigroup in which the idempotents form a subsemigroup. By \cite[Theorem VI.1.1]{key15}, the following conditions are equivalent for a regular semigroup $S$: (a) $S$ is orthodox; (b) $V_S(e)\subseteq E_S$ for all $e\in E_S$; (c) $V_S(b)V_S(a)\subseteq V_S(ab)$ for all $a, b\in S$. Thus if $S$ is orthodox and $x\!\perp\!y$ in $S$, then $x\!\in\! N_S$ if and only if $y\!\in\! N_S$, and $x^n\!\perp\!y^n$ for all $n\!\in\!\N$ (we denote by $\N$ the set of all positive integers), so that $o(x)=o(y)$. Moreover, if $S$ is any semigroup and $xy\!=\!x^2y^2$ for some $x,y\!\in\! S$,  it is clear that $xy\!=\!x^ny^n$ for all $n\!\in\!\N$. These simple facts will be used below without comment.

Let $S$ be an orthodox semigroup, and let $\cY=\{(x,y)\in S\times S:V(x)=V(y)\}$. Then $\cY$  is the smallest inverse semigroup congruence on $S$ \cite[Theorem VI.1.12]{key15}, so that $S/\cY$ is the maximum inverse semigroup homomorphic image of $S$. As noted, for instance, in \cite[p. 72]{key10}, it is easily seen (and well known) that $S$ is combinatorial if and only if $S/\cY$ is such.

In \cite{key7, key8}, a semigroup $S$ was termed {\em elementary} if $S=\langle A\cup B\rangle$ for two nonempty subsets $A$ and $B$ of $S$ such that $a\perp b$ for all $a\in A$ and $b\in B$, and an orthodox semigroup $S=\langle a, b\rangle$ with $a\perp b$ was referred to as ``an elementary orthodox semigroup on two mutually inverse generators''. This notion was introduced as a generalization of (and by analogy with) ``elementary generalized groups'' studied by Gluskin in \cite{key9}. In modern terminology, ``elementary generalized groups'' are known as ``monogenic inverse semigroups'' (by definition, an inverse semigroup is {\em monogenic} if it is generated by a pair of mutually inverse elements). Since the phrase ``an elementary orthodox semigroup on two mutually inverse generators'' is too cumbersome, we introduce a different term for such semigroups. Namely, we will say that $S$ is a {\em monogenic orthodox semigroup} (by analogy with the inverse semigroup case) if $S$ is an orthodox semigroup such that $S=\langle a, b\rangle$ for some $a, b\in S$ with $a\perp b$. Of course, it is important to require in this definition that $S$ be orthodox in addition to being generated by a pair of mutually inverse elements. By \cite[Exercise 2.3.5]{key5}, a semigroup given by the  presentation $\langle p, q\,|\,pqp=p, qpq=q\rangle$ is not regular; in fact, $\{p, q, pq, qp\}$ is its only regular $\cD$-class. Moreover, as shown by Clifford (see \cite[Example 1.6]{key7}), even if $a\perp b$ in a regular semigroup $S$, then $\langle a, b\rangle$ need not be a regular subsemigroup of $S$. On the other hand, an elementary subsemigroup of an orthodox semigroup is itself orthodox \cite[Proposition 1.1]{key7}, so if $S$ is orthodox and $a\perp b$ in $S$, then $\langle a, b\rangle$ is a monogenic orthodox semigroup. In what follows, the phrase ``let $S=\langle a, b\rangle$ be a monogenic orthodox semigroup'' will always mean that $S$ is an orthodox semigroup generated by a pair of mutually inverse elements $a, b\in S$. 

In view of the above, it is natural to ask: If $S=\langle x, y\rangle$ is a semigroup such that $x\perp y$, under what conditions is $S$ orthodox? The answer follows easily from \cite[Theorem 1.1]{key8}, stating that the free monogenic orthodox semigroup ${\rm FO}$ has the following presentation: ${\rm FO}=\langle p, q\,|\,p^n\perp q^n\text{ for all }n\in\N\rangle$; we will also denote this semigroup by ${\rm FO}(p,q)$ if we wish to name the generators explicitly. 
\bl\label{11} Let $S=\langle x, y\rangle$ be a semigroup such that $x\perp y$. Then
\vspace{0.02in}\\
\indent{\rm(i)} $S$ is orthodox if and only if $x^n\perp y^n$ for all $n\in\N$;
\vspace{0.02in}\\
\indent{\rm(ii)} $S$ is orthodox if $xy=x^2y^2$.
\el
{\bf Proof.} (i) If $x^n\perp y^n$ for all $n\in\N$, then $S$ is a homomorphic image of ${\rm FO}(x,y)$, and hence, by \cite[Lemma II.4.7]{key15}, $S$ is orthodox. The converse is obvious.
\vspace{0.01in}\\
\indent(ii) Suppose that $xy=x^2y^2$. Then for any $n\geq 1$,
$x^ny^nx^n=(x^ny^n)x^n=(xy)x^n=x^n$ and $y^nx^ny^n=y^n(x^ny^n)=y^n(xy)=y^n$.
Therefore, by (i), $S$ is orthodox.\epr
\vspace{0.05in}\\ 
\indent In \cite[Section II]{key8}, using the description of ${\rm FO}$ from \cite[Section I]{key8} and a number of results from \cite{key6} and a few other papers, Eberhart and Williams analyzed the lattice of congruences on ${\rm FO}$ and constructed several interesting monogenic orthodox semigroups. In particular, in \cite[Results 2.3, 2.4]{key8} they exhibited the following two infinitely presented bisimple monogenic orthodox semigroups: ${\rm FO}/\alpha'=\langle p, q\,|\,p^n\!\perp \!q^n\,(\forall\,n\!\in\!\N)\text{ and }pq\!=\!p^2q^2\rangle$ and ${\rm FO}/\sigma'=\langle p, q\,|\,p^n\!\perp\! q^n\,(\forall\,n\!\in\!\N)\text{ and }pq\!=\!p^2q^2, qp\!=\!q^2p^2\rangle$, where $\alpha'$ and $\sigma'$ are certain congruences on ${\rm FO}$ defined in \cite{key8} in several steps using the description of congruences on the free monogenic inverse semigroup and the results about the lattice of congruences on ${\rm FO}$ (the interested reader is referred to \cite{key8} for definitions of $\alpha'$ and $\sigma'$; they will not be needed in this paper). It is immediate from Lemma \ref{11} that ${\rm FO}/\alpha'$ and ${\rm FO}/\sigma'$ can be finitely presented:
\bc\label{12} The semigroups ${\rm FO}/\alpha'$ and ${\rm FO}/\sigma'$ given in \cite[Results 2.3 and 2.4]{key8}, respectively, have the following finite presentations: ${\rm FO}(p,q)/\alpha'=\langle\, p, q\;|\;p\perp q,\; pq=p^2q^2\,\rangle$ and ${\rm FO}(p,q)/\sigma'=\langle\, p, q\;|\;p\perp q,\; pq=p^2q^2,\; qp=q^2p^2\,\rangle$. 
\ec
It is stated without proof in \cite{key8} that each element of ${\rm FO}/\alpha'$ has a unique representation in the form $p^iq^mp^nq^j$, where $m, n\!\geq\! 0$, $i, j\!\in\!\{0,1\}$, $m\!\geq\! i$, and $n\!\geq\! j$ (the uniqueness assertion here is false since, for instance, $pq\!=\!p^1q^1p^0q^0\!=\!p^0q^0p^1q^1$). Note also that the eggbox picture of ${\rm FO}/\alpha'$ given in \cite[Figure 2]{key8} contains a number of misprints (among other things, $R_q$ is typed there twice). After correcting these inaccuracies, one could, in principle, obtain a description of all bisimple monogenic orthodox semigroups by using \cite[Results 2.3, 2.4]{key8} together with several other results about the lattice of congruences on ${\rm FO}$ stated in \cite[Section II]{key8} (mostly without proofs). Nevertheless \cite{key8} contains no such description, and at any rate it is certainly of interest to get it independently of all the theorems on the structure of the lattice of congruences on ${\rm FO}$ used in \cite{key8}. The main goal of this paper is to obtain such a description directly, using only some basic results about orthodox semigroups.

The following simple fact is probably well-known. We record it for convenience of reference and include its proof for completeness.
\bl\label{13} Let $S$ be an orthodox semigroup, $a$ an arbitrary element of $S$, and $b\in V(a)$. 
\vspace{0.02in}\\
\indent{\rm(i)} The following conditions are equivalent: {\rm(a)} $a^2$ is a left divisor of $a$; {\rm (b)} $a=a^2(b^2a)$; {\rm (c)} $ab=a^2b^2$; {\rm (d)} $b=(ba^2)b^2$; {\rm(e)} $b^2$ is a right divisor of $b$.
\vspace{0.02in}\\
\indent{\rm(ii)} If $ab=a^2b^2$ and $ba=b^2a^2$, then $H_a$ and $H_b$ are groups with identity elements $ab^2a$ and $ba^2b$, respectively. Conversely, if $H_a$ (and thus $H_b$) is a group, then $ab=a^2b^2$ and $ba=b^2a^2$.
\el 
{\bf Proof.} (i) Suppose $a^2$ is a left divisor of $a$. Then $a=a^2x$ for some $x\in S$, and hence $a^2(b^2a)=(a^2b^2a^2)x=a^2x=a$. Thus (a) implies (b). Since the converse is obvious, (a)$\Leftrightarrow$(b). By symmetry, (d)$\Leftrightarrow$(e), and the equivalences (b)$\Leftrightarrow$(c) and (c)$\Leftrightarrow$(d) hold trivially.
\vspace{0.02in}\\
\indent (ii) Suppose $ab=a^2b^2$ and $ba=b^2a^2$. Then $(ab^2a)a=a=a(ab^2a)$ whence $(a, ab^2a)\in\cR$ and $(a, ab^2a)\in\cL$, that is, $(a, ab^2a)\in\cH$. Since $(ab^2a)^2=a(b^2a^2b^2)a=ab^2a$, by Green's Theorem \cite[Theorem 2.16]{key5}, $H_a$ is a group and $ab^2a$ is its identity element. By symmetry, $H_b$ is a group with identity $ba^2b$. 

Conversely, assume that $H_a$ is a group. Since $b\perp a$ and $S$ is orthodox, $H_b$ is also a group and $b^2\perp a^2$. Thus, by \cite[Lemmas 2.12 and 2.15]{key5}, $a^2b^2=ab$ and $b^2a^2=ba$.  \epr
\bl\label{14}
Let $S$ be an orthodox semigroup and $a\perp b$ in $S$. The following conditions are equivalent: {\rm (a)} $ab=a^2b^2$ and $ba\not=b^2a^2$; {\rm (b)} $ab=a^2b^2$ and $a\in N_S$; {\rm (c)} $a^2$ is a left divisor of $a$ and $a\in N_S$; {\rm (d)} $a^2$ is a left divisor of $a$ and $ba\not=b^2a^2$. If one (and hence each) of these conditions holds, then $o(a)=o(b)=\infty$, $\langle a\rangle\cup\langle a\rangle b\subseteq R_a$, and $\langle b\rangle\cup a\langle b\rangle\subseteq L_b$.
\el 
{\bf Proof.} It is immediate from Lemma \ref{13} that conditions (a), (b), (c), and (d) are equivalent. Suppose (c) holds. Then $a=a^2x$ for some $x\in S$. If $o(a)<\infty$, there is $m\in\N$ such that $a^m\in E_S$ and, by putting $e=a^m$, we get $a=a^2x=a^3x^2=\cdots=a^mx^{m-1}=ex^{m-1}$ whence $a=ea=a^{m+1}$, which shows that $\langle a\rangle$ is a group; a contradiction. Therefore $o(a)\!=\!\infty\!=\!o(b)$. Note that $a\!=\!a^2x$ implies $a\cR a^2$, and so $a\cR a^n$ for any $n\in\N$. Moreover, $a^nb\cR a^n$ for all $n\in\N$ since $(a^nb)a=a^n$. Hence $\langle a\rangle\cup\langle a\rangle b\subseteq R_a$ and, by symmetry, $\langle b\rangle\cup a\langle b\rangle\subseteq L_b$.  \epr
\vspace{0.05in}\\
\indent The equivalence of conditions (a)--(d) of Lemma \ref{14} will be used below without mention.
\bl\label{15}
Let $S$ be an orthodox semigroup, $a\in N_S$, $b\in V(a)$, and $ab=a^2b^2$. Then
\vspace{0.02in}\\
\indent{\rm (i)} the subsemigroup $\langle ba^2, b^2a\rangle$ of $S$ is bicyclic and $ba$ is its identity element, that is, $\langle ba^2, b^2a\rangle=\cB(ba^2, b^2a)$; 
\vspace{0.02in}\\
\indent{\rm (ii)} $(ba^2)^k=ba^{k+1}$ and $(b^2a)^k=b^{k+1}a$ for all $k\geq 0$, so each element of $\langle ba^2, b^2a\rangle$ has a unique representation in the form $b^{m+1}a^{n+1}$ for some nonnegative integers $m$ and $n$;
\vspace{0.02in}\\
\indent{\rm (iii)} if $k$, $l$, $ m$, and $n$ are nonnegative integers such that $k+l\geq 1$ and $m+n\geq 1$, then $b^ka^l\cR b^ma^n$ if and only if $k\!=\!m$, and $b^ka^l\cL b^ma^n$ if and only if $l=n$.  
\el
{\bf Proof}. (i) This is immediate since $(ba^2)(b^2a)=ba$, $(b^2a)(ba^2)=b^2a^2\neq ba$, and $ba$ is a two-sided identity for $ba^2$ and $b^2a$. 
\vspace{0.02in}\\
\indent (ii) Take an arbitrary nonnegative integer $k$. If $k\geq 1$, then $(ba^2)^k=ba(aba)^{k-1}a=ba^{k+1}$ and $(b^2a)^k=b(bab)^{k-1}ba=b^{k+1}a$. Since $ba$ is the identity element of $\cB(ba^2, b^2a)$, we define $(ba^2)^0=(b^2a)^0=ba$, so the assertion holds for $k=0$ as well. Therefore each element of $\langle ba^2, b^2a\rangle$ can be uniquely written as $(b^2a)^m(ba^2)^n$, and thus as $(b^{m+1}a)(ba^{n+1})=b^{m+1}a^{n+1}$ for some nonnegative integers $m$ and $n$.  
\vspace{0.02in}\\
\indent (iii) Let $k$, $l$, $ m$, and $n$ be nonnegative integers such that $k+l\geq 1$ and $m+n\geq 1$. If $k, l, m, n\geq 1$, then $b^ka^l, b^ma^n\in\cB(ba^2, b^2a)$ and therefore $b^ka^l\cR b^ma^n$ if and only if $k=m$, and $b^ka^l\cL b^ma^n$ if and only if $l=n$.

Suppose $l=0$ (and hence $k\geq 1$). Let $b^k\cR b^ma^n$. If $m,n\geq 1$, using the fact that $b^k\cR b^ka$, we get $b^ka\cR b^ma^n$, so that $k=m$. If $n=0$ (and thus $m\geq 1$), then $b^k\cR b^m$ and so, assuming without loss of generality that $k\leq m$, we obtain $b=b(a^kb^k)\cR b(a^kb^k)b^{m-k}=b^{m-k+1}$ whence $m=k$ in view of Green's Theorem \cite[Theorem 2.16]{key5}, Lemma \ref{14}, and the fact that $b\in N_S$. Note that we cannot have $m=0$ and $n\geq 1$, for otherwise $b^k\cR a^n$ which implies $b^{k+1}\cR ba^n$ and thus, as shown above, $k=0$; a contradiction. We have proved that $b^k\cR b^ma^n$ implies $k=m$, and if $m\geq 1$, it is easily seen that $b^m\cR b^ma^n$ for all $n\geq 0$, so the converse also holds. Now let $b^k\cL b^ma^n$. If $m\geq 1$, then $b^ka\cL b^ma^{n+1}$ and since $b^ka,b^ma^{n+1}\in\cB(ba^2, b^2a)$, we have $n=0$. Here again we cannot have $m=0$ and $n\geq 1$, for otherwise $b^k\cL a^n$ and, using the fact that $a^n\cL ba^n$, we obtain $b^k\cL ba^n$ so that, as has just been shown, $n=0$; a contradiction. Thus $b^k\cL b^ma^n$ implies $n=0$, and the converse follows from Lemma \ref{14}. By symmetry, the result also holds in the case when $l\geq 1$ and $k=0$. The proof is complete.\epr   
\med
\section{Bisimple monogenic orthodox semigroups with nongroup generators}
\med
\indent Let $S$ be an orthodox semigroup, $a\in N_S$, $b\in V(a)$, and $ab=a^2b^2$. Clearly, every $x\!\in\!S$ can be written in the form $x\!=\!(a^{k_1}b^{l_1})\cdots(a^{k_r}b^{l_r})$ for some $r\!\in\!\N$ and nonnegative integers $k_i, l_i\;(i\!=\!1,\ldots,r)$ such that $k_i+l_i\!\geq\!1$ for all $1\!\leq\!i\!\leq\!r$. We call $a^{k_1}b^{l_1}$, \ldots, $a^{k_r}b^{l_r}$ the {\em syllables} of $(a^{k_1}b^{l_1})\cdots(a^{k_r}b^{l_r})$ and refer to $(a^{k_1}b^{l_1})\cdots(a^{k_r}b^{l_r})$ as an {\em $r$-syllable word} in $a, b$. If $m, n\in\Z$ and $i, j\in\{0, 1\}$ are such that either (I) $m>i$ and $n=j=0$, or (II) $m=i=0$ and $n\geq 1$, or (III) $m>i$ and $n>j$, we say that $a^ib^ma^nb^j$ is an {\em abridged word in $a,b$} (in this order!) of {\em type I}, {\em II}, or {\em III}, respectively (or simply an {\em abridged word in $a,b$} when there is no need to indicate its type), and if $x=a^ib^ma^nb^j$, we will also call $a^ib^ma^nb^j$ an {\em abridged form} of $x$.
\vspace{0.03in}\\ 
\indent{\bf Remark 1.} Let $S=\langle a, b\rangle$ be a monogenic orthodox semigroup with $ab=a^2b^2$ and $ba\neq b^2a^2$. As we will see shortly, each element of $S$ can be written in an abridged form although, in general, not uniquely. Later, by imposing certain additional conditions on $S$, we will associate with each $x\!\in\! S$ a unique abridged word for which we reserve the term ``reduced form'' of $x$. Observe also that definitions of type I and type II abridged words in $a,b$ are not entirely symmetric -- in a type I word $a^ib^m$ we must have $m>i$ but in a type II word $a^nb^j$ the equality $n=j=1$ is allowed. As will be seen later, this ensures that $ab\in S$ is always represented by a {\em unique} abridged word in $a,b$ (namely, by a type II word $a^0b^0a^1b^1$). Thus the dual of an abridged word $b^1a^1$ of type II in $b,a$ is an abridged word $a^1b^1$ of type II in $a,b$; this is the only exception of the following rule -- the dual of an abridged word in $b,a$ of type I [type II, type III] is an abridged word in $a,b$ of type II [type I, type III].
\bl\label{20} 
Let $S$ be an orthodox semigroup, $a\in N_S$, $b\in V(a)$, and $ab=a^2b^2$. Then
\vspace{0.02in}\\
\indent{\rm(i)} each $1$-syllable word in $a,b$ can be written as an abridged word in $a,b$ of type I or II;
\vspace{0.02in}\\
\indent{\rm (ii)} the product of two abridged words in $a,b$ of type I [type II] can be written as an abridged word in $a,b$ of type I [type II];
\vspace{0.02in}\\
\indent{\rm (iii)} if $x=a^ib^m$ and $y=a^nb^j$ are abridged words in $a,b$ of types I and II, respectively, then $xy$ equals an abridged word in $a,b$ of type $I$ or type III, while $yx$ can be written as a $1$-syllable word in $a,b$ and thus as an abridged word in $a,b$ of type I or type II.
\el
{\bf Proof.} (i) Let $x\!=\!a^kb^l$ for some $k,l\geq 0$ with $k+l\geq 1$. If $k=l$, then $x$ equals $ab$, an abridged word in $a,b$ of type II. If $k\!=\!0$ or $l\!=\!0$, then $x$ equals $b^l$ or $a^k$, an abridged word in $a,b$ of type I or II, respectively. If $l\!>\!k\!\geq\!1$, then $x\!=\!a^kb^kb^{l-k}\!=\!ab^{l-k+1}$, an abridged word in $a,b$ of type I, and if $k\!>\!l\!\geq\!1$, then $x\!=\!a^{k-l}a^lb^l\!=\!a^{k-l+1}b$, an abridged word in $a,b$ of type II.
\vspace{0.02in}\\
\indent (ii) If $a^ib^m$ and $a^{i'}b^{m'}$ are abridged words in $a,b$ of type I, then $a^ib^ma^{i'}b^{m'}=a^ib^{m+m'-i'}$, which is an abridged word in $a,b$ of type I. Similarly, if $a^nb^j$ and $a^{n'}b^{j'}$ are abridged words in $a,b$ of type II, then $a^nb^ja^{n'}b^{j'}=a^{n+n'-j}b^{j'}$, which is an abridged word in $a,b$ of type II.
\vspace{0.02in}\\
\indent (iii) Let $x\!=\!a^ib^m$ and $y\!=\!a^nb^j$ be abridged words in $a,b$ of types I and II, respectively. If $n\!=\!j\!=\!1$, then $xy\!=\!a^ib^mab\!=\!a^ib^m$, and if $n\!>\!j$, then $xy\!=\!a^ib^ma^nb^j$, which, by definition, is an abridged word in $a,b$ of type III. Finally, $yx\!=\!a^nb^m$ if $i\!=\!j$, and $yx\!=\!a^{n+1-j}b^{m+1-i}$ if $i\!\neq\! j$, so that, by (i), $yx$ can be written as an abridged word in $a,b$ of type I or type II. \epr 
\bl\label{21} 
Let $S\!=\!\langle a, b\rangle$ be a monogenic orthodox semigroup such that $a\!\in\!N_S$ and $ab\!=\!a^2b^2$. Then for each $x\!\in\! S$ there is an abridged word $a^ib^ma^nb^j$ in $a, b$ such that $x\!=\!a^ib^ma^nb^j$.
\el
{\bf Proof.} Let $x\!\in\! S$, and let $r$ be the smallest positive integer such that $x\!=\!(a^{k_1}b^{l_1}\!)\cdots(a^{k_r}b^{l_r}\!)$ for some $r$-syllable word $(a^{k_1}b^{l_1})\cdots(a^{k_r}b^{l_r})$ in $a,b$. In view of Lemma \ref{20}(i), we can assume that $a^{k_i}b^{l_i}$ is an abridged word in $a,b$ of type I or II for each $i=1,\ldots,r$. Since $x$ cannot be represented by a $q$-syllable word for $q<r$, by Lemma \ref{20}(ii) adjacent syllables $a^{k_j}b^{l_j}$ and $a^{k_{j+1}}b^{l_{j+1}}\,(j=1,\ldots,r-1)$ are of different types, and by Lemma \ref{20}(iii) a syllable written as an abridged word in $a,b$ of type II cannot precede a syllable which is an abridged word in $a,b$ of type I. Therefore either $r=1$ and $a^{k_1}b^{l_1}$ is an abridged word in $a,b$ of type I or II, or $r=2$ and $a^{k_1}b^{l_1}a^{k_2}b^{l_2}$ is an abridged word in $a,b$ of type III. 
\epr
\vspace{0.03in}\\ 
\indent{\bf Remark 2.} By \cite[Lemma 1.3]{key8}, each element of ${\rm FO}(a,b)$ can be written uniquely as a word $x^{i_1}_1\cdots x^{i_n}_n$ in $a,b$, where $i_k\!>\!\min\{i_{k-1},i_{k+1}\}$ for $k\!=\!2,\ldots,n-1$ (and, of course, $x_j\!\neq\! x_{j+1}$ for $j=1,\ldots,n-1$). Lemma \ref{21} can be proved using that result but when this alternative proof is rigorously written, it is neither shorter nor better than the proof given above.   
\bl\label{22} Let $S=\langle a, b\rangle$ be a monogenic orthodox semigroup such that $ab=a^2b^2$ and $ba\neq b^2a^2$. Then  
\vspace{0.02in}\\
\indent{\rm(i)} if $a^ib^ma^nb^j$ is an abridged word in $a,b$ of type III (that is, $i,j\in\{0,1\}$, $m>i$, $n>j$), then $a^ib^m\cR a^ib^ma^nb^j$ and $a^ib^ma^nb^j\cL a^nb^j$;  
\vspace{0.02in}\\
\indent{\rm(ii)} $S$ is bisimple and combinatorial.  
\el
{\bf Proof.} (i) Suppose $a^ib^ma^nb^j$ is an abridged word in $a,b$ of type III. It is easily seen that $a^ib^m=a^ib^ma^nb^ja^jb^n$, which implies $a^ib^m\cR a^ib^ma^nb^j$. Similarly, since $a^nb^j=a^mb^ia^ib^ma^nb^j$, we have $a^ib^ma^nb^j\cL a^nb^j$.  
\vspace{0.03in}\\
\indent{\rm(ii)} Let $s\!\in\! S$. By Lemma \ref{21}, $s=a^ib^ma^nb^j$ where $a^ib^ma^nb^j$ is an abridged word in $a,b$. If $a^ib^ma^nb^j$ is of type I, $n\!=\!j\!=\!0$ and $m\!>\!i\!\in\!\{0,1\}$ whence $s\!=\!a^ib^m\cL b\cD a$ by Lemma \ref{14}. If $a^ib^ma^nb^j$ is of type II, $m\!=\!i\!=\!0$ and $n\!\geq\! 1$, so that $s\!=\!a^nb^j\cR a$ again by Lemma \ref{14}. Finally, if $a^ib^ma^nb^j$ is of type III, by part (i) and Lemma \ref{14}, $a^ib^ma^nb^j\cL a^nb^j\cR a$. Thus $S=D_a$. 

Let $x=a\cY$ and $y=b\cY$. Then $S/\cY=\langle x,y\rangle$, $x\perp y$ in $S/\cY$, and $xy=x^2y^2$, which implies $(xy)y=(xy)(yx)y=(yx)(xy)y=y(x^2y^2)=y$, and, by symmetry, $x(xy)=x$. Therefore $xy$ is the identity of $S/\cY$. If $xy=yx$, then $(ab)\cY=(ba)\cY$, so that $ab\perp ba$ in $S$ whence $ba=b(ab)(ba)(ab)a=b^2a^2$; a contradiction. Thus $yx\neq xy$ and, by \cite[Lemma 1.31]{key5}, $S/\cY=\cB(x,y)$. It follows that $S$ is combinatorial. \epr
\vspace{0.1in}\\
\indent {\bf Remark 3.} In the proof of Lemma \ref{22}(ii), it was shown that $S/\cY=\cB(x,y)$, and we could have deduced from this that $S$ is bisimple by applying, for instance, \cite[Theorem 10]{key13}. However, a simple direct proof of that fact given above seems to be preferable. 
\vspace{0.07in}\\
\indent Let $S=\langle a, b\rangle$ be a monogenic orthodox semigroup such that $ab=a^2b^2$ and $ba\not=b^2a^2$. If $n\in\N$, it is clear that $(ab)a^n=a^n$ and $b^n(ab)=b^n$ but for the products $a^n(ab)$ and $(ab)b^n$ in $S$ there are the following possibilities:
\vspace{0.06in}\\
\indent 1) $a^n(ab)\not=a^n$ and $(ab)b^n\not=b^n$ for all $n\in\N$, in which case we denote $S$ by $\cO_{(\infty, \infty)}(a,b)$;  
\vspace{0.04in}\\
\indent 2) $(ab)b^m\!\not=\!b^m$ for all $m\!\in\!\N$ but $a^k(ab)=a^k$ for some $k\in\N$, and letting $n$ be the smallest of such integers $k$, we denote $S$ by $\cO_{(n, \infty)}(a,b)$;
\vspace{0.04in}\\
\indent 3) $a^n(ab)\not=a^n$ for all $n\in\N$ but $(ab)b^l=b^l$ for some $l\in\N$, and with $m$ standing for the smallest of such integers $l$, we denote $S$ by $\cO_{(\infty, m)}(a,b)$; 
\vspace{0.04in}\\
\indent 4) $a^n(ab)=a^n$ and $(ab)b^m=b^m$ for some $n,m\in\N$, and letting $n$ and $m$ be the smallest integers with these properties, we denote $S$ by $\cO_{(n, m)}(a,b)$ (of course, $\cO_{(1, 1)}(a,b)$ is just another notation for the bicyclic semigroup $\cB(a,b)$).
\vspace{0.1in}\\
\indent {\bf Remark 4.} In view of Remark 1, the dual of $\cO_{(\infty, \infty)}(a,b)$ is $\cO_{(\infty, \infty)}(b,a)$, the dual of $\cO_{(n, \infty)}(a,b)$ is $\cO_{(\infty, n)}(b,a)$, and the dual of $\cO_{(n, m)}(a,b)$ is $\cO_{(m, n)}(b,a)$. Thus all results about $\cO_{(\infty, m)}(a,b)$ are obtained automatically from the corresponding results about $\cO_{(m, \infty)}(b,a)$.   
\bp\label{23}
Let $S\!=\!\langle a, b\rangle$ be a monogenic orthodox semigroup such that $a\!\in\!N_S$ and $ab\!=\!a^2b^2$. Then every nontrivial relation that holds in $S$ (other than $ab=a^2b^2$) is equivalent either to $a^{n+1}b=a^n$, or to $ab^{m+1}=b^m$, or to $a^{n+1}b=a^n$ and $ab^{m+1}=b^m$ for some $m, n\in\N$. 
\ep
{\bf Proof.} By Lemma \ref{14}, $o(a)=o(b)=\infty$. Suppose that $S$ satisfies a nontrivial relation $w(a,b)=w'(a,b)$, which is not equivalent to $ab=a^2b^2$. In view of Lemma \ref{21}, we may assume that $w(a,b)$ and $w'(a,b)$ are abridged words in $a,b$.
\vspace{0.02in}\\
\indent{\bf Case 1:} $w(a,b)$ and $w'(a,b)$ are both of type I.
\vspace{0.01in}\\
\indent In this case, $w(a,b)=a^ib^m$ and $w'(a,b)=a^{i'}b^{m'}$ with $i,i'\in\{0,1\}$, $m>i$, and $m'>i'$. If $i\!=\!i'$, then $b^m\!=\!b^{m'}$ whence $m\!=\!m'$, so that $w(a,b)\!=\!w'(a,b)$ is a trivial relation; a contradiction. Therefore $i\!\neq\! i'$. Assume without loss of generality that $i\!=\!0$ and $i'\!=\!1$. Then $b^m\!=\!ab^{m'}$ whence $b^{m+1}\!=\!b^{m'}$ and so $m'\!=\!m+1$. Thus the given relation is $ab^{m+1}\!=\!b^m$.
\vspace{0.02in}\\
\indent{\bf Case 2:} $w(a,b)$ is of type I and $w'(a,b)$ of type II.
\vspace{0.01in}\\
\indent Here we have $a^ib^m=a^{n'}b^{j'}$ with $i,j'\in\{0,1\}$, $m>i$, and $n'\geq 1$. Thus $b^{m-i+1}=ba^{n'}b^{j'}$, which implies $b^{(m-i)+(n'-j')+1}=b^{m-i+1}b^{n'-j'}=ba^{n'}b^{n'}=b$. Hence $(m-i)+(n'-j')=0$, contrary to the conditions $m-i>0$ and $n'-j'\geq 0$. Therefore this case cannot happen.
\vspace{0.02in}\\
\indent{\bf Case 3:} $w(a,b)$ is of type I and $w'(a,b)$ of type III.
\vspace{0.01in}\\
\indent In this case, $a^ib^m\!=\!a^{i'}b^{m'}a^{n'}b^{j'}$ where $i,i',j'\!\in\!\{0,1\}$, $m\!>\!i$, $m'\!>\!i'$, and $n'\!>\!j'$. Then $b(a^ib^m)a\!=\!b(a^{i'}b^{m'}a^{n'}b^{j'})a$, that is, $b^{m-i+1}a\!=\!b^{m'-i'+1}a^{n'-j'+1}$. By Lemma \ref{15}, $m-i\!=\!m'-i'$ and $n'-j'\!=\!0$, contradicting the fact that $n'\!>\!j'$. Hence this case is also impossible.   
\vspace{0.02in}\\
\indent{\bf Case 4:} $w(a,b)$ and $w'(a,b)$ are both of type II.
\vspace{0.01in}\\ 
\indent Here $a^nb^j=a^{n'}b^{j'}$ with $j,j'\in\{0,1\}$, $n\geq 1$, and $n'\geq 1$, and hence $a^nb^ja=a^{n'}b^{j'}a$, that is, $a^{n-j+1}=a^{n'-j'+1}$, which implies $n-j=n'-j'$. Since the relation $a^nb^j=a^{n'}b^{j'}$ is nontrivial, we have $j\neq j'$. Without loss of generality assume that $j=0$ and $j'=1$. Then $a^n=a^{n'}b$ whence $a^{n+1}=a^{n'}$ and therefore $n'=n+1$. Thus the given relation is $a^{n+1}b=a^n$.
\vspace{0.02in}\\
\indent{\bf Case 5:} $w(a,b)$ is of type II and $w'(a,b)$ of type III.
\vspace{0.01in}\\
\indent In this case, $a^nb^j=a^{i'}b^{m'}a^{n'}b^{j'}$ with $j,i',j'\in\{0,1\}$, $n\geq 1$, $m'>i'$, and $n'>j'$. Then $b(a^nb^j)a=b(a^{i'}b^{m'}a^{n'}b^{j'})a$, that is, $ba^{n-j+1}=b^{m'-i'+1}a^{n'-j'+1}$. By Lemma \ref{15}, $n-j=n'-j'$ and $m'-i'\!=\!0$, but the latter contradicts the fact that $m'\!>\!i'$. Thus this case is impossible.
\vspace{0.02in}\\
\indent{\bf Case 6:} $w(a,b)$ and $w'(a,b)$ are both of type III. 
\vspace{0.01in}\\
\indent The given relation is $a^ib^ma^nb^j=a^{i'}b^{m'}a^{n'}b^{j'}$ with $i,j,i',j'\in\{0,1\}$, $m>i$, $n>j$, $m'>i'$, and $n'>j'$. Then $b(a^ib^ma^nb^j)a=b(a^{i'}b^{m'}a^{n'}b^{j'})a$, that is, $b^{m-i+1}a^{n-j+1}=b^{m'-i'+1}a^{n'-j'+1}$ and so, by Lemma \ref{15}, $m-i=m'-i'$ and $n-j=n'-j'$. Since the given relation is nontrivial, we cannot have both $i=i'$ and $j=j'$.
\vspace{0.02in}\\
\indent 6a) Suppose that $i=i'$ and $j\not=j'$. Assuming without loss of generality that $j=0$ and $j'=1$, we have the relation $a^ib^ma^n=a^ib^{m'}a^{n'}b$. Hence $b(a^ib^ma^n)a=b(a^ib^{m'}a^{n'}b)a$, that is, $b^{m-i+1}a^{n+1}=b^{m'-i+1}a^{n'}$. By Lemma \ref{15}, $m'=m$ and $n'=n+1$, so $a^ib^ma^n=a^ib^{m'}a^{n'}b$ is, in fact, $a^ib^ma^n=a^ib^ma^{n+1}b$. Then $a^n=a^mb^ma^n=a^{m-i}a^ib^ma^n=a^{m-i}a^ib^ma^{n+1}b=a^{n+1}b$. Since $a^n\!=\!a^{n+1}b$ implies $a^ib^ma^n\!=\!a^ib^ma^{n+1}b$, the given relation is equivalent to $a^n\!=\!a^{n+1}b$.
\vspace{0.02in}\\
\indent 6b) By symmetry with 6a), if $i\neq i'$ and $j=j'$, the relation is equivalent to $b^m=ab^{m+1}$.
\vspace{0.02in}\\
\indent 6c) Suppose that $i\neq i'$ and $j\not=j'$. Assume without loss of generality that $i=0$ and $i'=1$. Suppose, first, that $j=0$ and $j'=1$. Since $b^ma^n=ab^{m'}a^{n'}b$ implies $b^{m+1}a^{n+1}=b^{m'}a^{n'}$ and so, by Lemma \ref{15}, $m'=m+1$ and $n'=n+1$, the given relation is $b^ma^n=ab^{m+1}a^{n+1}b$. Hence $b^m=b^ma^nb^n=ab^{m+1}a^{n+1}b^{n+1}=ab^{m+1}$ and $a^n=a^mb^ma^n=a^{m+1}b^{m+1}a^{n+1}b=a^{n+1}b$. On the other hand, $b^m=ab^{m+1}$ and $a^n=a^{n+1}b$ together clearly imply $b^ma^n=ab^{m+1}a^{n+1}b$. Therefore the given relation is equivalent to two relations, $b^m=ab^{m+1}$ and $a^n=a^{n+1}b$. Assume now that $j=1$ and $j'=0$. Since $b^ma^nb=ab^{m'}a^{n'}$ implies $b^{m+1}a^n=b^{m'}a^{n'+1}$, by Lemma \ref{15}, $m'=m+1$ and $n=n'+1$, so the given relation is $b^ma^nb=ab^{m+1}a^{n-1}$. Then $b^m=b^ma^nb^n=ab^{m+1}a^{n-1}b^{n-1}=ab^{m+1}$ and $a^{n-1}=a^{m+1}b^{m+1}a^{n-1}=a^mb^ma^nb=a^nb$. Since $b^m=ab^{m+1}$ and $a^{n-1}=a^nb$ together clearly imply $b^ma^nb=ab^{m+1}a^{n-1}$, the given relation is equivalent to two relations, $b^m=ab^{m+1}$ and $a^{n-1}=a^nb$. This completes the proof. \epr 
\vspace{0.05in}\\
\indent It is asserted (without proof) in \cite{key8} that \cite[Figure 2]{key8} gives an eggbox picture of ${\rm FO}(a, b)/\alpha'$ but, as pointed out in the paragraph following Corollary \ref{12}, \cite[Figure 2]{key8} contains quite a few misprints. Here is its corrected version (with $p$ and $q$ replaced by $a$ and $b$, respectively): 
\vspace{0.01in}\\
\[
\renewcommand{\arraystretch}{1.2}
\begin{tabular}{c|c|c|c|c|c|c|c|c} 
$\vdots\;\vdots\;\vdots$&$\vdots$&$\vdots$&$\vdots$&$\vdots$&$\vdots$&$\vdots$&$\vdots$&$\vdots\;\vdots\;\vdots$\\ \hline
$\cdots$&$ab^4a^4b$&$ab^4a^3b$&$ab^4a^2b$&$ab^4$&$ab^4a$&$ab^4a^2$&$ab^4a^3$&$\cdots$\\ \hline
$\cdots$&$ab^3a^4b$&$ab^3a^3b$&$ab^3a^2b$&$ab^3$&$ab^3a$&$ab^3a^2$&$ab^3a^3$&$\cdots$\\ \hline
$\cdots$&$ab^2a^4b$&$ab^2a^3b$&$ab^2a^2b$&$ab^2$&$ab^2a$&$ab^2a^2$&$ab^2a^3$&$\cdots$\\ \hline
$\cdots$&$a^4b$&$a^3b$ & $a^2b$  & $ab$    & $a$     & $a^2$   & $a^3$      &$\cdots$\\ \hline
$\cdots$&$ba^4b$&$ba^3b$ & $ba^2b$ & $b$ & $ba$ & $ba^2$ & $ba^3$           &$\cdots$\\ \hline
$\cdots$&$b^2a^4b$&$b^2a^3b$&$b^2a^2b$&$b^2$ & $b^2a$  & $b^2a^2$& $b^2a^3$ &$\cdots$\\ \hline
$\cdots$&$b^3a^4b$&$b^3a^3b$&$b^3a^2b$&$b^3$ & $b^3a$  & $b^3a^2$& $b^3a^3$ &$\cdots$\\ \hline
$\vdots\;\vdots\;\vdots$&$\vdots$&$\vdots$&$\vdots$&$\vdots$&$\vdots$&$\vdots$&$\vdots$&$\vdots\;\vdots\;\vdots$\\ 
\end{tabular} 
\]
\vspace{-0.1in}\\
\begin{center}Figure 1\vspace{0.02in}\\
\end{center}
It follows from Proposition \ref{23} that $\cO_{(\infty, \infty)}(a,b)=\langle\, a, b\;|\;a\perp b,\, ab=a^2b^2\,\rangle$. Therefore, by Corollary \ref{12}, $\cO_{(\infty, \infty)}(a,b)={\rm FO}(a, b)/\alpha'$ which means (in view of the above mentioned assertion in \cite{key8} about the eggbox picture of ${\rm FO}(a, b)/\alpha'$) that Figure 1 exhibits the eggbox picture of  $\cO_{(\infty, \infty)}(a,b)$. We will include a proof of this fact both for the sake of completeness and because it is independent of the results about the lattice of congruences on ${\rm FO}$ contained in \cite[Section II]{key8}.  
\bl\label{24} The semigroup $\cO_{(\infty, \infty)}(a,b)$ is a combinatorial bisimple monogenic orthodox semigroup whose eggbox picture is given in Figure 1.
\el
{\bf Proof.} Let $S=\cO_{(\infty, \infty)}(a,b)$. According to Lemma \ref{22}, $S$ is bisimple and combinatorial. The latter implies that $R_a\cap L_b=\{ab\}$. By Lemma \ref{14}, $\langle a\rangle\cup\langle a\rangle b\subseteq R_a$ and $\langle b\rangle\cup a\langle b\rangle\subseteq L_b$. It is clear that $\langle a\rangle\cup\langle a\rangle b$ coincides with the set of all abridged words in $a,b$ of type II, and from the definition of $\cO_{(\infty, \infty)}(a,b)$ and the proof of Proposition \ref{23} (see Case 4), it is immediate that $a^nb^j=a^{n'}b^{j'}$ if and only if $n=n'$ and $j=j'$ for all $a^nb^j, a^{n'}b^{j'}\in\langle a\rangle\cup\langle a\rangle b$. It follows also from the definition of $\cO_{(\infty, \infty)}(a,b)$ and the proof of Proposition \ref{23} that if $i,i'\in\{0,1\}$ and $k,k'\in\N$, then $a^ib^k=a^{i'}b^{k'}$ if and only if $i=i'$ and $k=k'$, and it is obvious that the set of all abridged words in $a,b$ of type I coincides with $(\langle b\rangle\cup a\langle b\rangle)\setminus\{ab\}$. Suppose that $a^ib^ma^nb^j$ is an abridged word of type I or III such that $a^ib^ma^nb^j\in R_a$. By Lemma \ref{22}(i), $a^ib^m\cR a^ib^ma^nb^j$ and hence $a^ib^m\in R_a\cap L_b=\{ab\}$; a contradiction. Therefore $R_a=\langle a\rangle\cup\langle a\rangle b$ and, similarly, $L_b=\langle b\rangle\cup a\langle b\rangle$, which also shows that $S\setminus(R_a\cup L_b)$ coincides with the set of all abridged words in $a,b$ of type III. Finally, if $a^ib^ma^nb^j$ is any abridged word of type III, Lemma \ref{22}(i) and the above remarks imply that $\{a^ib^ma^nb^j\}=R_{a^ib^m}\cap L_{a^nb^j}$. It follows that Figure 1 gives the eggbox picture of $\cO_{(\infty, \infty)}(a,b)$. The proof is complete.\epr
\vspace{0.02in}\\
\indent{\bf Remark 5.} It is immediate from Lemma \ref{24} that for each $x\!\in\!\cO_{(\infty, \infty)}(a,b)$ there is a unique abridged word $a^ib^ma^nb^j$ in $a,b$ such that $x\!=\!a^ib^ma^nb^j$, and we will say that $a^ib^ma^nb^j$ is a {\em reduced word representing} $x$ (or the {\em reduced form} of $x$). It is also easily seen that $\cO_{(\infty, \infty)}(a,b)$ is the disjoint union of four bicyclic semigroups and two infinite cyclic semigroups: $\cO_{(\infty, \infty)}(a,b)\!=\!\cB(a^2b,ab^2)\cup\cB(ab^2a^2,ab^3a)\cup\cB(ba^3b,b^2a^2b)\cup\cB(ba^2,b^2a)\cup\langle a\rangle\cup\langle b\rangle$.
\vspace{0.02in}\\
\indent Let $n\in\N$. By definition of $\cO_{(n, \infty)}(a,b)$, we have $(ab)b^m\neq b^m$ for all $m\in\N$, $a^n(ab)=a^n$, and $a^l(ab)\neq a^l$ for all $1\leq l< n$. Thus $\cO_{(n, \infty)}(a,b)=\langle a,b\,|\,a\perp b,\, ab=a^2b^2,\,a^{n+1}b=a^n\,\rangle$, which shows that $\cO_{(n, \infty)}(a,b)$ is a homomorphic image of $\cO_{(\infty, \infty)}(a,b)$, and the eggbox picture of $\cO_{(n, \infty)}(a,b)$ is easily obtained from that of $\cO_{(\infty, \infty)}(a,b)$; it is given in Figure 2.
\vspace{0.01in}\\
\[
\renewcommand{\arraystretch}{1.2}
\begin{tabular}{|c|c|c|c|c|c|c|c|c} 
$\vdots$&$\vdots$&$\vdots\;\vdots\;\vdots$&$\vdots$&$\vdots$&$\vdots$&$\vdots$&$\vdots$&$\vdots\;\vdots\;\vdots$\\ \hline
$ab^4a^nb$&$ab^4a^{n-1}b$&$\cdots$&$ab^4a^2b$&$ab^4$&$ab^4a$&$ab^4a^2$&$ab^4a^3$&$\cdots$\\ \hline
$ab^3a^nb$&$ab^3a^{n-1}b$&$\cdots$&$ab^3a^2b$&$ab^3$&$ab^3a$&$ab^3a^2$&$ab^3a^3$&$\cdots$\\ \hline
$ab^2a^nb$&$ab^2a^{n-1}b$&$\cdots$&$ab^2a^2b$&$ab^2$&$ab^2a$&$ab^2a^2$&$ab^2a^3$&$\cdots$\\ \hline
$a^nb$&$a^{n-1}b$&$\cdots $&$a^2b$ & $ab$ & $a$  & $a^2$   & $a^3$    & $\cdots$\\ \hline
$ba^nb$&$ba^{n-1}b$&$\cdots $ & $ba^2b$ & $b$ & $ba$ & $ba^2$ & $ba^3$           &$\cdots$\\ \hline
$b^2a^nb$&$b^2a^{n-1}b$&$\cdots$&$b^2a^2b$&$b^2$ & $b^2a$  & $b^2a^2$& $b^2a^3$ &$\cdots$\\ \hline
$b^3a^nb$&$b^3a^{n-1}b$&$\cdots$&$b^3a^2b$&$b^3$ & $b^3a$  & $b^3a^2$& $b^3a^3$ &$\cdots$\\ \hline
$\vdots$&$\vdots$&$\vdots\;\vdots\;\vdots$&$\vdots$&$\vdots$&$\vdots$&$\vdots$&$\vdots$&$\vdots\;\vdots\;\vdots$\\ 
\end{tabular} 
\]
\vspace{-0.1in}\\
\begin{center}Figure 2\vspace{0.02in}\\
\end{center}
Note that for each $x\in\cO_{(n, \infty)}(a,b)$ there is a unique abridged word $a^ib^ma^lb^j$ with $l\leq n$ such that $x=a^ib^ma^lb^j$; we will call it the {\em reduced form} of $x$. By duality, if $m\in\N$, we obtain the eggbox picture of $\cO_{(\infty, m)}(a,b)$ and observe that each $x\in\cO_{(\infty, m)}(a,b)$ is represented by a unique abridged word $a^ib^ka^nb^j$ with $k\leq m$, called the {\em reduced form} of $x$. 
\vspace{0.02in}\\
\indent Now take any $m,n\in\N$. By definition of $\cO_{(n, m)}(a,b)$, we have $a^n(ab)=a^n$ and $(ab)b^m=b^m$  but $a^l(ab)\neq a^l$ for all $1\leq l< n$ and $(ab)b^k\neq b^k$ for all $1\leq k<m$. Therefore $\cO_{(n, m)}(a,b)$ is given by the presentation: $\cO_{(n, m)}(a,b)=\langle a,b\,|\,a\perp b,\, ab=a^2b^2,\,a^{n+1}b=a^n,\,ab^{m+1}=b^m\,\rangle$, and so, in particular, $\cO_{(n, m)}(a,b)$ is a homomorphic image of $\cO_{(n, \infty)}(a,b)$. This enables us to obtain the eggbox picture of $\cO_{(n, m)}(a,b)$ from that of $\cO_{(n, \infty)}(a,b)$; it is shown in Figure 3. It is clear that for every $x\in\cO_{(n, m)}(a,b)$ there is a unique abridged word $a^ib^ka^lb^j$ such that $k\leq m$, $l\leq n$, and $x=a^ib^ka^lb^j$; we will say that $a^ib^ka^lb^j$ is the {\em reduced form} of $x$.
\vspace{0.05in}\\
\[
\renewcommand{\arraystretch}{1.2}
\begin{tabular}{|c|c|c|c|c|c|c|c} 
\hline
$ab^ma^nb$&$ab^ma^{n-1}b$&$\cdots$&$ab^ma^2b$&$ab^m$&$ab^ma$&$ab^ma^2$&$\cdots$\\ \hline
$ab^{m-1}a^nb$&$ab^{m-1}a^{n-1}b$&$\cdots$&$ab^{m-1}a^2b$&$ab^{m-1}$&$ab^{m-1}a$&$ab^{m-1}a^2$&$\cdots$\\ \hline
$\vdots$&$\vdots$&$\vdots\;\vdots\;\vdots$&$\vdots$&$\vdots$&$\vdots$&$\vdots$&$\vdots\;\vdots\;\vdots$\\ \hline
$ab^3a^nb$&$ab^3a^{n-1}b$&$\cdots$&$ab^3a^2b$&$ab^3$&$ab^3a$&$ab^3a^2$&$\cdots$\\ \hline
$ab^2a^nb$&$ab^2a^{n-1}b$&$\cdots$&$ab^2a^2b$&$ab^2$&$ab^2a$&$ab^2a^2$&$\cdots$\\ \hline
$a^nb$&$a^{n-1}b$&$\cdots$ & $a^2b$  & $ab$    & $a$     & $a^2$   & $\cdots$\\ \hline
$ba^nb$&$ba^{n-1}b$&$\cdots$ & $ba^2b$ & $b$ & $ba$ & $ba^2$ & $\cdots$\\ \hline
$b^2a^nb$&$b^2a^{n-1}b$&$\cdots$&$b^2a^2b$&$b^2$ & $b^2a$  & $b^2a^2$& $\cdots$\\ \hline
$\vdots$&$\vdots$&$\vdots\;\vdots\;\vdots$&$\vdots$&$\vdots$&$\vdots$&$\vdots$&$\vdots\;\vdots\;\vdots$\\ 
\end{tabular} 
\]
\vspace{0.02in}\\
\begin{center}Figure 3\vspace{0.1in}\\\end{center}

For convenience of reference, we summarize the above observations in the following lemma.   
\bl\label{25} Let $m,n\in\N$. The semigroups $\cO_{(n, \infty)}(a,b)$ and $\cO_{(n, m)}(a,b)$ are combinatorial bisimple monogenic orthodox semigroups whose eggbox pictures are given in Figures 2 and 3, respectively, and the semigroup $\cO_{(\infty, m)}(a,b)$ is the dual of $\cO_{(m, \infty)}(b,a)$.
\el 

Recall that a band $E$ is said to be {\em uniform} if $eEe\cong fEf$ for all $e,f\in E$. As shown by Hall \cite[Main Theorem]{key12}, a band $E$ is the band of a bisimple orthodox semigroup if and only if $E$ is uniform. In particular, the easy part of the cited theorem guarantees that the band of idempotents of any bisimple orthodox semigroup is uniform (see \cite[Result 7]{key12} or \cite[Proposition VI.3.1]{key15}). Therefore, for any $m,n\in\N$, the bands of idempotents of semigroups $\cO_{(\infty, \infty)}(a,b)$, $\cO_{(n, \infty)}(a,b)$, $\cO_{(\infty, m)}(a,b)$, and $\cO_{(n, m)}(a,b)$ are uniform, and using the diagrams of these bands we can explicitly illustrate this fact. 

In view of duality, it is sufficient to consider the bands of idempotents of the semigroups $\cO_{(\infty, \infty)}(a,b)$, $\cO_{(n, \infty)}(a,b)$, and $\cO_{(n, m)}(a,b)$ with $m\geq n$; they are shown in parts (a), (b), and (c), respectively, of Figure 4 with the bold line segments representing the covering relation of the natural partial order and the thin line segments indicating the $\cR$- and $\cL$-relations on each of these bands. (The diagrams of $\cO_{(n, \infty)}(a,b)$ and $\cO_{(n, m)}(a,b)$ are drawn under the assumptions that $n>1$ and $m>n>1$, respectively; modifications for $n=1$ and for $m=n>1$ or $m>n=1$ are obvious.) It is easily seen that if $E$ is any of the bands shown in Figure 4 and $e\in E$, then $eEe$ is isomorphic to the chain of idempotents of the bicyclic semigroup (that is, to the chain $e_0>e_1>e_2>\cdots$), and hence $eEe\cong\! fEf$ for all $e,f\in E$.      
\vspace{0.05in}\\
\begin{picture}(420,250)

\put(-5,138){\line(5,4){44}}
\put(-5,84){\line(5,4){44}}
\put(75,138){\line(5,4){44}}
\put(75,84){\line(5,4){44}}

\put(40,224){\line(0,-1){150}}
\put(39.7,224){\line(0,-1){150}}
\put(39.3,224){\line(0,-1){150}}
\put(40.3,224){\line(0,-1){150}}
\put(40.7,224){\line(0,-1){150}}

\put(-5,138){\line(0,-1){85}}
\put(-4.7,138){\line(0,-1){85}}
\put(-4.3,138){\line(0,-1){85}}
\put(-5.3,138){\line(0,-1){85}}
\put(-5.7,138){\line(0,-1){85}}

\put(75,138){\line(0,-1){85}}
\put(74.7,138){\line(0,-1){85}}
\put(74.3,138){\line(0,-1){85}}
\put(75.3,138){\line(0,-1){85}}
\put(75.7,138){\line(0,-1){85}}

\put(120,174){\line(0,-1){85}}
\put(119.7,174){\line(0,-1){85}}
\put(119.3,174){\line(0,-1){85}}
\put(120.3,174){\line(0,-1){85}}
\put(120.7,174){\line(0,-1){85}}

\put(-5,138){\line(1,0){80}}
\put(-5,84){\line(1,0){80}}

\put(82, 176){\tiny{$\cR$}}
\put(50, 140){\tiny{$\cR$}}
\put(101, 154){\tiny{$\cL$}}
\put(21, 154){\tiny{$\cL$}}

\put(50, 86){\tiny{$\cR$}}
\put(82, 122){\tiny{$\cR$}}
\put(101, 100){\tiny{$\cL$}}
\put(21, 100){\tiny{$\cL$}}

\put(40,174){\line(1,0){80}}
\put(40,120){\line(1,0){80}}

\put(40,224){\circle*{5}}

\put(40,174){\circle*{5}} 
\put(120,174){\circle*{5}} 

\put(40,120){\circle*{5}}
\put(120,120){\circle*{5}}

\put(-5,138){\circle*{5}}
\put(75,138){\circle*{5}}

\put(-5,84){\circle*{5}}
\put(75,84){\circle*{5}}

\put(27,225){\tiny{$ab$}} 

\put(10,174){\tiny{$ab^2a^2b$}}
\put(114,179){\tiny{$ab^2a$}}

\put(10,117){\tiny{$ab^3a^3b$}}
\put(123,123){\tiny{$ab^3a^2$}}

\put(-2,129.5){\tiny{$ba^2b$}}
\put(80,136){\tiny{$ba$}}
\put(-2,75.5){\tiny{$b^2a^3b$}}
\put(80,82){\tiny{$b^2a^2$}}

\put(39,55){$\vdots$}
\put(-6,35){$\vdots$}
\put(74,35){$\vdots$}
\put(119,70){$\vdots$}

\put(45,-38){$\text{(a)}$}


\put(165,158){\line(5,4){44}}
\put(165,84){\line(5,4){44}}
\put(245,158){\line(5,4){44}}
\put(245,84){\line(5,4){44}}
\put(245,36){\line(5,4){44}}
\put(245,6){\line(5,4){44}}

\put(210,234){\line(0,-1){55}}
\put(209.7,234){\line(0,-1){55}}
\put(209.3,234){\line(0,-1){55}}
\put(210.3,234){\line(0,-1){55}}
\put(210.7,234){\line(0,-1){55}}

\put(165,158){\line(0,-1){15}}
\put(164.7,158){\line(0,-1){15}}
\put(164.3,158){\line(0,-1){15}}
\put(165.3,158){\line(0,-1){15}}
\put(165.7,158){\line(0,-1){15}}

\put(245,158){\line(0,-1){15}}
\put(244.7,158){\line(0,-1){15}}
\put(244.3,158){\line(0,-1){15}}
\put(245.3,158){\line(0,-1){15}}
\put(245.7,158){\line(0,-1){15}}

\put(290,194){\line(0,-1){15}}
\put(289.7,194){\line(0,-1){15}}
\put(289.3,194){\line(0,-1){15}}
\put(290.3,194){\line(0,-1){15}}
\put(290.7,194){\line(0,-1){15}}

\put(165,84){\line(0,1){10}}
\put(164.7,84){\line(0,1){10}}
\put(164.3,84){\line(0,1){10}}
\put(165.3,84){\line(0,1){10}}
\put(165.7,84){\line(0,1){10}}

\put(210,120){\line(5,-3){80}}
\put(210.2,120){\line(5,-3){80}}
\put(210.5,120){\line(5,-3){80}}
\put(210.8,120){\line(5,-3){80}}
\put(211,120){\line(5,-3){80}}
\put(209,120){\line(5,-3){80}}
\put(209.2,120){\line(5,-3){80}}
\put(209.5,120){\line(5,-3){80}}
\put(209.8,120){\line(5,-3){80}}
\put(208.8,120){\line(5,-3){80}}
\put(208.5,120){\line(5,-3){80}}
\put(208.3,120){\line(5,-3){80}}
\put(208,120){\line(5,-3){80}}

\put(165,84){\line(5,-3){80}}
\put(165.2,84){\line(5,-3){80}}
\put(165.5,84){\line(5,-3){80}}
\put(165.8,84){\line(5,-3){80}}
\put(166,84){\line(5,-3){80}}
\put(164,84){\line(5,-3){80}}
\put(164.2,84){\line(5,-3){80}}
\put(164.5,84){\line(5,-3){80}}
\put(164.8,84){\line(5,-3){80}}
\put(163.8,84){\line(5,-3){80}}
\put(163.5,84){\line(5,-3){80}}
\put(163.3,84){\line(5,-3){80}}
\put(163,84){\line(5,-3){80}}

\put(245,84){\line(0,1){10}}
\put(244.7,84){\line(0,1){10}}
\put(244.3,84){\line(0,1){10}}
\put(245.3,84){\line(0,1){10}}
\put(245.7,84){\line(0,1){10}}

\put(210,120){\line(0,1){10}}
\put(209.7,120){\line(0,1){10}}
\put(209.3,120){\line(0,1){10}}
\put(210.3,120){\line(0,1){10}}
\put(210.7,120){\line(0,1){10}}

\put(290,120){\line(0,1){10}}
\put(289.7,120){\line(0,1){10}}
\put(289.3,120){\line(0,1){10}}
\put(290.3,120){\line(0,1){10}}
\put(290.7,120){\line(0,1){10}}

\put(245,84){\line(0,-1){86}}
\put(244.7,84){\line(0,-1){86}}
\put(244.3,84){\line(0,-1){86}}
\put(245.3,84){\line(0,-1){86}}
\put(245.7,84){\line(0,-1){86}}

\put(290,118){\line(0,-1){84}}
\put(289.7,118){\line(0,-1){84}}
\put(289.3,118){\line(0,-1){84}}
\put(290.3,118){\line(0,-1){84}}
\put(290.7,118){\line(0,-1){84}}

\put(165,158){\line(1,0){80}}
\put(165,84){\line(1,0){80}}

\put(252, 196){\tiny{$\cR$}}
\put(216, 160){\tiny{$\cR$}}
\put(273, 175){\tiny{$\cL$}}
\put(193, 175){\tiny{$\cL$}}

\put(198, 87){\tiny{$\cR$}}
\put(234, 123){\tiny{$\cR$}}
\put(272, 100){\tiny{$\cL$}}
\put(192, 100){\tiny{$\cL$}}

\put(267, 60){\tiny{$\cL$}}
\put(267, 30){\tiny{$\cL$}}

\put(210,194){\line(1,0){80}}
\put(210,120){\line(1,0){80}}

\put(210,234){\circle*{5}}

\put(210,194){\circle*{5}} 
\put(290,194){\circle*{5}} 

\put(210,120){\circle*{5}}
\put(290,120){\circle*{5}}

\put(165,158){\circle*{5}}
\put(245,158){\circle*{5}}

\put(165,84){\circle*{5}}
\put(245,84){\circle*{5}}

\put(245,36){\circle*{5}}
\put(245,6){\circle*{5}}

\put(290,72){\circle*{5}}
\put(290,42){\circle*{5}}

\put(197,235){\tiny{$ab$}} 

\put(178,194){\tiny{$ab^2a^2b$}}
\put(280,199){\tiny{$ab^2a$}}

\put(178,120){\tiny{$ab^na^nb$}}
\put(294,120){\tiny{$ab^n\!a^{\!n-\!1}$}}

\put(294,72){\tiny{$ab^{n+1}\!a^n$}}

\put(294,42){\tiny{$ab^{n+2}a^{n+1}$}}

\put(170,148){\tiny{$ba^2b$}}
\put(250,155){\tiny{$ba$}}
\put(140,73){\tiny{$b^{n-\!1}a^nb$}}
\put(209,75){\tiny{$b^{n-\!1}\!a^{n-\!1}$}}
\put(222,30){\tiny{$b^na^n$}}
\put(210,3){\tiny{$b^{n\!+\!1}a^{n\!+\!1}$}}

\put(209,163){$\vdots$}
\put(209,141){$\vdots$}
\put(164,128){$\vdots$}
\put(164,106){$\vdots$}
\put(244,128){$\vdots$}
\put(244,106){$\vdots$}
\put(289,163){$\vdots$}
\put(289,141){$\vdots$}

\put(244,-15){$\vdots$}
\put(289,22){$\vdots$}

\put(220,-38){$\text{(b)}$}


\put(335,183){\line(5,4){44}}
\put(335,109){\line(5,4){44}}
\put(415,183){\line(5,4){44}}
\put(415,109){\line(5,4){44}}

\put(415,76){\line(5,4){44}}

\put(415,40){\line(5,4){44}}

\put(380,244){\line(0,-1){40}}
\put(379.7,244){\line(0,-1){40}}
\put(379.3,244){\line(0,-1){40}}
\put(380.3,244){\line(0,-1){40}}
\put(380.7,244){\line(0,-1){40}}

\put(335,183){\line(0,-1){15}}
\put(334.7,183){\line(0,-1){15}}
\put(334.3,183){\line(0,-1){15}}
\put(335.3,183){\line(0,-1){15}}
\put(335.7,183){\line(0,-1){15}}

\put(415,183){\line(0,-1){15}}
\put(414.7,183){\line(0,-1){15}}
\put(414.3,183){\line(0,-1){15}}
\put(415.3,183){\line(0,-1){15}}
\put(415.7,183){\line(0,-1){15}}

\put(460,219){\line(0,-1){15}}
\put(459.7,219){\line(0,-1){15}}
\put(459.3,219){\line(0,-1){15}}
\put(460.3,219){\line(0,-1){15}}
\put(460.7,219){\line(0,-1){15}}

\put(335,109){\line(0,1){10}}
\put(334.7,109){\line(0,1){10}}
\put(334.3,109){\line(0,1){10}}
\put(335.3,109){\line(0,1){10}}
\put(335.7,109){\line(0,1){10}}

\put(380,145){\line(0,1){10}}
\put(379.7,145){\line(0,1){10}}
\put(379.3,145){\line(0,1){10}}
\put(380.3,145){\line(0,1){10}}
\put(380.7,145){\line(0,1){10}}

\put(460,145){\line(0,1){10}}
\put(459.7,145){\line(0,1){10}}
\put(459.3,145){\line(0,1){10}}
\put(460.3,145){\line(0,1){10}}
\put(460.7,145){\line(0,1){10}}

\put(415,120){\line(0,-1){53}}
\put(414.7,120){\line(0,-1){53}}
\put(414.3,120){\line(0,-1){53}}
\put(415.3,120){\line(0,-1){53}}
\put(415.7,120){\line(0,-1){53}}

\put(460,145){\line(0,-1){43}}
\put(459.7,145){\line(0,-1){43}}
\put(459.3,145){\line(0,-1){43}}
\put(460.3,145){\line(0,-1){43}}
\put(460.7,145){\line(0,-1){43}}

\put(335,183){\line(1,0){80}}
\put(335,109){\line(1,0){80}}

\put(335,109){\line(5,-2){80}}
\put(335.3,109){\line(5,-2){80}}
\put(335.5,109){\line(5,-2){80}}
\put(335.8,109){\line(5,-2){80}}
\put(336,109){\line(5,-2){80}}
\put(334,109){\line(5,-2){80}}
\put(334.2,109){\line(5,-2){80}}
\put(334.5,109){\line(5,-2){80}}
\put(334.7,109){\line(5,-2){80}}
\put(333,109){\line(5,-2){80}}
\put(333.2,109){\line(5,-2){80}}
\put(333.5,109){\line(5,-2){80}}
\put(333.7,109){\line(5,-2){80}}
\put(332.6,109){\line(5,-2){80}}
\put(332.8,109){\line(5,-2){80}}

\put(380,145){\line(5,-2){80}}
\put(380.3,145){\line(5,-2){80}}
\put(380.5,145){\line(5,-2){80}}
\put(380.8,145){\line(5,-2){80}}
\put(379,145){\line(5,-2){80}}
\put(379.2,145){\line(5,-2){80}}
\put(379.5,145){\line(5,-2){80}}
\put(379.7,145){\line(5,-2){80}}
\put(378,145){\line(5,-2){80}}
\put(378.2,145){\line(5,-2){80}}
\put(378.5,145){\line(5,-2){80}}
\put(378.7,145){\line(5,-2){80}}
\put(377.6,145){\line(5,-2){80}}
\put(377.8,145){\line(5,-2){80}}

\put(415,14){\line(3,4){45}}
\put(415.2,14){\line(3,4){45}}
\put(415.5,14){\line(3,4){45}}
\put(415.8,14){\line(3,4){45}}
\put(414.2,14){\line(3,4){45}}
\put(414.5,14){\line(3,4){45}}
\put(414.8,14){\line(3,4){45}}
\put(414,14){\line(3,4){45}}

\put(415,48){\line(0,-1){60}}
\put(414.7,48){\line(0,-1){60}}
\put(414.3,48){\line(0,-1){60}}
\put(415.3,48){\line(0,-1){60}}
\put(415.7,48){\line(0,-1){60}}

\put(460,83){\line(0,-1){6}}
\put(459.7,83){\line(0,-1){6}}
\put(459.3,83){\line(0,-1){6}}
\put(460.3,83){\line(0,-1){6}}
\put(460.7,83){\line(0,-1){6}}

\put(418, 221){\tiny{$\cR$}}
\put(387, 185){\tiny{$\cR$}}
\put(442, 200){\tiny{$\cL$}}
\put(362, 200){\tiny{$\cL$}}

\put(368, 111){\tiny{$\cR$}}
\put(406, 147){\tiny{$\cR$}}
\put(444, 127){\tiny{$\cL$}}
\put(364, 127){\tiny{$\cL$}}

\put(427, 92){\tiny{$\cL$}}

\put(427, 56){\tiny{$\cL$}}

\put(380,219){\line(1,0){80}}
\put(380,145){\line(1,0){80}}

\put(380,244){\circle*{5}}

\put(380,219){\circle*{5}} 
\put(460,219){\circle*{5}} 

\put(380,145){\circle*{5}}
\put(460,145){\circle*{5}}

\put(335,183){\circle*{5}}
\put(415,183){\circle*{5}}

\put(335,109){\circle*{5}}
\put(415,109){\circle*{5}}

\put(415,77){\circle*{5}}
\put(460,112){\circle*{5}}

\put(460,75){\circle*{5}}

\put(415,40){\circle*{5}}

\put(415,15){\circle*{5}}

\put(415,-5){\circle*{5}}

\put(367,244){\tiny{$ab$}} 

\put(348,219){\tiny{$ab^2a^2b$}}
\put(450,224){\tiny{$ab^2a$}}

\put(348,145){\tiny{$ab^na^nb$}}
\put(429,148){\tiny{$ab^n\!a^{n-\!1}$}}

\put(426,107){\tiny{$ab^{n+1}\!a^n$}}

\put(340,173){\tiny{$ba^2b$}}
\put(420,179){\tiny{$ba$}}
\put(312,99){\tiny{$b^{n-\!1}a^nb$}}
\put(380,100){\tiny{$b^{n-\!1}\!a^{n-\!1}$}}

\put(393,71){\tiny{$b^na^n$}}

\put(425,74){\tiny{$ab^{m}\!a^{m-\!1}$}}

\put(373,38){\tiny{$b^{m-\!1}a^{m-\!1}$}}

\put(421,14){\tiny{$b^{m}a^{m}$}}

\put(421,-6){\tiny{$b^{m+1}a^{m+1}$}}

\put(379,188){$\vdots$}
\put(379,170){$\vdots$}

\put(334,152){$\vdots$}
\put(334,134){$\vdots$}

\put(414,152){$\vdots$}
\put(414,134){$\vdots$}

\put(459,188){$\vdots$}
\put(459,170){$\vdots$}

\put(414,54){$\vdots$}
\put(459,88){$\vdots$}

\put(414,-24){$\vdots$}

\put(388,-38){$\text{(c)}$}

\put(203,-58){$\text{Figure }4$}
\end{picture}
\vspace{1in}\\
\indent We have shown that if a monogenic orthodox semigroup $S=\langle a,b\rangle$ such that $a,b\in N_S$ satisfies $ab=a^2b^2$, then $S$ coincides with one of the bisimple combinatorial semigroups $\cO_{(\infty, \infty)}(a,b)$, $\cO_{(n, \infty)}(a,b)$, $\cO_{(\infty, m)}(a,b)$, or $\cO_{(n, m)}(a,b)$ for some $m,n\in\N$, described by Lemmas \ref{24} and \ref{25}. We will prove now that there are no other bisimple orthodox semigroups generated by two mutually inverse nongroup elements and thus obtain a complete description of all monogenic orthodox semigroups with nongroup generators. 

By \cite[Lemma 2.7]{key10} (deduced in \cite{key10} as a corollary to \cite[Section II]{key8}), if $S$ is an orthodox semigroup, $a\!\in\! N_S$, and $b\!\in\! V(a)$, then either $\{a, b, ab, ba\}$ is the top $\cD$-class of $\langle a, b\rangle$ (which means that $\langle a,b\rangle\setminus\{a, b, ab, ba\}$ is an ideal of $\langle a, b\rangle$), or $o(a)=\infty=o(b)$. As mentioned prior to the statement of \cite[Lemma 2.7]{key10}, it could have been proved independently of \cite{key8}. Now we will give such a direct proof establishing, in fact, a sharper result. 
\bl\label{26}{\rm (An extension of \cite[Lemma 2.7]{key10})} Let $S$ be an orthodox semigroup, $a$ an arbitrary nongroup element of $S$, and $b\in V(a)$. Then either $\{a, b, ab, ba\}$ is the top $\cD$-class of $\langle a, b\rangle$, or (perhaps after interchanging $a$ and $b$) $ab=a^2b^2$ and $ba\not=b^2a^2$, and hence $o(a)=o(b)=\infty$, $\langle a\rangle\cup\langle a\rangle b\subseteq R_a$, and $\langle b\rangle\cup a\langle b\rangle\subseteq L_b$.
\el
{\bf Proof}. Denote $A=\langle a, b\rangle$, and put, as always, $I^A(a)=J^A(a)\setminus J^A_a$ (the superscript $A$ is used to distinguish subsets of $A$ from the corresponding subsets of $S$). Suppose $a^2\in I^A(a)$. Since $I^A(a)$ is an ideal of $A$, we conclude that
\( J^A_a=J^A(a)\setminus I^A(a)=\{a, b, ab, ba\}\subseteq D^A_a\subseteq J^A_a.\)
Thus $J^A_a=D^A_a=\{a, b, ab, ba\}$. Moreover, $A\setminus \{a, b, ab, ba\}$ coincides with $I^A(a)$ and hence is an ideal of $A$. Therefore $\{a, b, ab, ba\}$ is the top $\cD$-class of $\langle a, b\rangle$. For the rest of the proof, we will assume that $a^2\not\in I^A(a)$, that is, $(a, a^2)\in\cJ^A$.

Suppose $a^2$ is neither a left nor a right divisor of $a$. By assumption, $a\in Aa^2A$. Therefore we can choose shortest possible words $u$ and $v$ in $a,b$ such that $a=ua^2v$. Since $b=b(ua^2v)b$ and since, by Lemma \ref{13}, $b^2$ is neither a left nor a right divisor of $b$, the first letter of $u$ as well as the last letter of $v$ is $a$. In fact, in view of our assumption about $a^2$, we must have $u=abu'$ and $v=v'ba$ where $u'$ and $v'$ are nonempty words in $a, b$. The first letter of $u'$ cannot be $a$, for otherwise $u'=au''$ for some nonempty word $u''$ in $a, b$, and hence $a=(au'')a^2v$ with $au''$ being shorter than $u$, a contradiction. Similarly, the last letter of $v'$ cannot be $a$. Thus $b$ is the first letter of $u'$ and the last letter of $v'$. However, $b=bu'a^2v'b$ and so $b^2$ is both a left and a right divisor of $b$, again a contradiction. 

Therefore $a^2$ is either a left or a right divisor of $a$ and, by Lemma \ref{13}, $b^2$ is, respectively, a right or a left divisor of $b$. Thus, interchanging $a$ and $b$ if necessary, we may assume that $a^2$ is a left divisor of $a$. Then, according to Lemma \ref{13}, $ab=a^2b^2$ and $ba\not=b^2a^2$ whence, by Lemma \ref{14}, $o(a)=o(b)=\infty$, $\langle a\rangle\cup\langle a\rangle b\subseteq R_a$, and $\langle b\rangle\cup a\langle b\rangle\subseteq L_b$.\epr
\vspace{0.04in}\\
\indent The author gratefully acknowledges that he learned the argument used in the second paragraph of the above proof from Mark Sapir (in a private conversation many years ago).
\bc\label{27} Let $S=\langle a,b\rangle$ be a bisimple monogenic orthodox semigroup such that $a$ (and therefore $b$) is a nongroup element of $S$. Then (perhaps after interchanging $a$ and $b$) we have $ab=a^2b^2$ and $ba\neq b^2a^2$.
\ec
The principal results of this section are combined in the following theorem.
\bt\label{28} Let $S=\langle a, b\rangle$ be a monogenic orthodox semigroup such that $a$ (and hence $b$) is a nongroup element of $S$. Then $S$ is bisimple if and only if $S$ (or the dual of $S$) coincides with one of the semigroups $\cO_{(\infty, \infty)}(a,b)$, $\cO_{(n, \infty)}(a,b)$, $\cO_{(\infty, m)}(a,b)$, or $\cO_{(n, m)}(a,b)$ for some $m,n\in\N$.
\et     
\sm
\section{Monogenic orthodox semigroups with group generators}
\sm
Let $S=\langle a, b\rangle$ be a monogenic orthodox semigroup such that $a$ (and hence $b$) is a group element of $S$. {\em To avoid repetition, this assumption will be retained through the rest of this section.} By Lemma \ref{13}(ii), $ab=a^2b^2$ and $ba=b^2a^2$, and the identity elements of the groups $H_a$ and $H_b$ are $ab^2a$ and $ba^2b$, respectively.  Since $ab=a^2b^2$ and $ba=b^2a^2$, it is easily seen that each element of $S$ can be written (in general, not uniquely) in the form $b^ia^nb^j$ or $a^ib^na^j$ for some $i,j\in\{0,1\}$ and $n\in\N$. Denote $e_a=ab^2a$ and $e_b=ba^2b$. Then $e_ae_b=ab\in R_a\cap L_b$ and $e_be_a=ba\in L_a\cap R_b$, and hence $E_S$ is a rectangular band with at most four elements: $e_a$, $e_b$, $ab$, and $ba$. It follows that $S$ is isomorphic to the direct product of $H_a$ and $E_S$ (see, for example, \cite[Exercise III.12]{key15}) so, in particular, $S$ is bisimple. There are four possibilities:
\vspace{0.03in}\\
\indent 1) $a(ab)\neq a$ and $(ab)b\neq b$, which happens precisely when the idempotents $e_a$, $e_b$, $ab$, and $ba$ are pairwise distinct, and thus $E_S=\{e_a, ab, ba, e_b\}$ is a $2\times 2$ rectangular band;
\vspace{0.02in}\\
\indent 2) $a(ab)\neq a$ and $(ab)b=b$, which is equivalent to the condition $e_a=ba\neq ab=e_b$, and hence $E_S=\{e_a, e_b\}$ is a $2$-element right zero semigroup;
\vspace{0.02in}\\
\indent 3) $a(ab)=a$ and $(ab)b\neq b$, which is equivalent to $e_a=ab\neq ba=e_b$, so that $E_S=\{e_a, e_b\}$ is a $2$-element left zero semigroup;
\vspace{0.02in}\\
\indent 4) $a(ab)=a$ and $(ab)b=b$, that is, $e_a=ab=ba=e_b$, and thus $E_S$ is a singleton.
\vspace{0.02in}\\
The eggbox pictures of $S$ in cases 1 -- 4 are shown, respectively, in parts (a) -- (d) of Figure 5.
\vspace{0.02in}\\ 
\[
\renewcommand{\arraystretch}{1.2}
\begin{tabular}{|c|c|} 
\hline
$H_a$&$H_{ab}$\\ \hline
$H_{ba}$&$H_b$\\ \hline
\end{tabular} 
\hspace{0.7in}
\begin{tabular}{|c|c|} 
\hline
$H_a$&$H_b$\\ \hline
\end{tabular}
\hspace{0.7in}
\begin{tabular}{|c|} 
\hline
$H_a$\\ \hline
$H_b$\\ \hline
\end{tabular}
\hspace{0.7in}
\begin{tabular}{|c|} 
\hline
$H_a$\\ \hline
\end{tabular}  
\]
\vspace{0.005in}\\
\hspace*{1.45in}(a)\hspace{1.2in}(b)\hspace{1in(c)}\hspace{0.85in}(d)
\vspace{-0.1in}\\
\begin{center}Figure 5\vspace{0.02in}\\
\end{center}  

{\bf Case 1:} $a(ab)\neq a$ and $(ab)b\neq b$. 
\vspace{0.02in}\\
\indent Suppose that $o(a)=o(b)=\infty$. Since $ab^3a$ is the inverse of $a$ in the group $H_a$ and $(ab^3a)^n\!=\!ab^{n+2}a$ for all $n\!\in\!\N$, we have $H_a\!\supseteq\!\langle a, ab^3a\rangle\!=\!\{\ldots,ab^4a,ab^3a,ab^2a,a,a^2,\ldots\}\!\cong\!\Z$. Dually, $H_b\supseteq\langle b,ba^3b\rangle=\{\ldots,ba^4b,ba^3b,ba^2b,b,b^2,\ldots\}\cong\Z$. Note that $a^2b$ is the inverse of $ab^2$ in the group $H_{ab}$, and $(a^2b)^n=a^{n+1}b$ and $(ab^2)^n=ab^{n+1}$ for all $n\in\N$. It follows that $H_{ab}\supseteq\{\ldots,ab^3,ab^2,ab,a^2b,a^3b,\ldots\}=\langle a^2b,ab^2\rangle\cong\Z$. By symmetry, we also have $H_{ba}\supseteq\{\ldots,ba^3,ba^2,ba,b^2a,b^3a,\ldots\}=\langle b^2a,ba^2\rangle\cong\Z$. Since each element of $S$ has the form $b^ia^nb^j$ or $a^ib^na^j$ for some $i,j\in\{0,1\}$ and $n\in\N$, we conclude that $H_a=\langle a, ab^3a\rangle$, $H_{ab}=\langle a^2b,ab^2\rangle$, $H_{ba}=\langle b^2a,ba^2\rangle$, and $H_b=\langle b,ba^3b\rangle$. Note that in this case $S$ coincides with the semigroup ${\rm FO}(a,b)/\sigma'$ of \cite[Result 2.4]{key8}.

Now let $o(a)\!=\!o(b)\!=\!m\!\in\!\N$. Then it is clear that $H_a, H_b, H_{ab}$, and $H_{ba}$ are cyclic groups of order $m$. More specifically, $H_a\!=\!\{ab^2a,a,\ldots,a^{m-1}\}\!=\!\langle a\rangle$, $H_b=\{ba^2b,b,\ldots,b^{m-1}\}=\langle b\rangle$, $H_{ab}=\{ab,a^2b,\ldots,a^{m-1}b\}=\langle a^2b\rangle$, and $H_{ba}=\{ba,b^2a,\ldots,b^{m-1}a\}=\langle b^2a\rangle$.
\vspace{0.02in}\\ 
\indent{\bf Case 2:} $a(ab)\neq a$ and $(ab)b=b$. 
\vspace{0.02in}\\
\indent Here, as observed above, $ab^2a\!=\!ba\!\neq\! ab\!=\!ba^2b$, so $E_S=\{ba,ab\}$ is a $2$-element right zero semigroup, and $S\!\cong\! H_a\times E_S$, that is, $S$ is a right group. Note that $b^na\!=\!ab^{n+1}a$ and $a^n\!=\!ba^{n+1}$ for all $n\!\in\!\N$. If $o(a)\!=\!o(b)\!=\!\infty$, then $H_a\!=\!H_{ba}\!=\!\langle a,b^2a\rangle\!=\!\{\ldots,b^3a,b^2a,ba,a,a^2,\ldots\}\!\cong\!\Z$ and $H_b\!=\!H_{ab}\!=\!\langle b,a^2b\rangle\!=\!\{\ldots,a^3b,a^2b,ab,b,b^2,\ldots\}\!\cong\!\Z$. Now assume that $o(a)\!=\!o(b)\!=\!m\!\in\!\N$. Then $H_a=H_{ba}\cong\Z_m\cong H_{ab}=H_b$ or, in more detail, $H_a=H_{ba}=\{ba,a,\ldots,a^{m-1}\}$ and $H_{ab}=H_b=\{ab,b,\ldots,b^{m-1}\}$.
\vspace{0.02in}\\ 
\indent{\bf Case 3:} $a(ab)=a$ and $(ab)b\neq b$.  
\vspace{0.02in}\\
\indent Here the description of the structure of $S$ is obtained from that of Case 2 by duality.
\vspace{0.02in}\\ 
\indent{\bf Case 4:} $a(ab)=a$ and $(ab)b=b$.
\vspace{0.02in}\\  
\indent In this case, $E_S$ is a singleton and $S=H_a$ is a cyclic group. 
\vspace{0.1in}\\ 
\indent The results of this section are summarized in the following theorem. 
\bt Let $S=\langle a,b\rangle$ be a monogenic orthodox semigroup such that $a$ (and hence $b$) is a group element of $S$. 

{\rm(1)} Let $a(ab)\neq a$ and $(ab)b\neq b$. Then $E_S=\{ab^2a,ab,ba,ba^2b\}$ is a $2\times 2$ rectangular band and $S\cong H_a\times E_S$ is a bisimple semigroup whose eggbox picture is shown in Figure 5(a). If $o(a)=\infty$, then $H_a=\langle a, ab^3a\rangle$, $H_{ab}=\langle a^2b,ab^2\rangle$, $H_{ba}=\langle b^2a,ba^2\rangle$, and $H_b=\langle b,ba^3b\rangle$ are infinite cyclic groups, and if $o(a)=m$ for some $m\in\N$, then the $\cH$-classes of $S$ are finite cyclic groups of order $m$: $H_a=\{ab^2a,a,\ldots,a^{m-1}\}=\langle a\rangle$, $H_b=\{ba^2b,b,\ldots,b^{m-1}\}=\langle b\rangle$, $H_{ab}=\{ab,a^2b,\ldots,a^{m-1}b\}=\langle a^2b\rangle$, and $H_{ba}=\{ba,b^2a,\ldots,b^{m-1}a\}=\langle b^2a\rangle$.

{\rm(2)} Let $a(ab)\neq a$ and $(ab)b=b$. Then $E_S=\{ba,ab\}$ is a $2$-element right zero semigroup, and $S\cong H_a\times E_S$ is a right group whose eggbox picture is shown in Figure 5(b). If $o(a)=\infty$, then $H_a=\langle a,b^2a\rangle$ and $H_b=\langle b,a^2b\rangle$ are infinite cyclic groups, and if $o(a)=m\in\N$, then $H_a=\{ba,a,\ldots,a^{m-1}\}=\langle a\rangle\cong\Z_m$ and $H_b=\{ab,b,\ldots,b^{m-1}\}=\langle b\rangle\cong\Z_m$.

{\rm(3)} If $a(ab)=a$ and $(ab)b\neq b$, then the structure of the dual of $S$ is described in {\rm (2)}.

{\rm(4)} If $a(ab)=a$ and $(ab)b=b$, then $S$ is a cyclic group.
\et
{\bf Remark 6.} The author is grateful to a reader of the preprint for proposing to use in the proof of Lemma \ref{22}(ii) the fact that if $S=\langle a,b\rangle$ is a monogenic orthodox semigroup with $ab=a^2b^2$ and $ba\neq b^2a^2$, then $S/\cY$ is bicyclic, and for suggesting that combinatorial bisimple monogenic orthodox semigroups can be described by investigating the idempotent-pure congruences on the semigroup $\langle a,b\,|\,a^i\perp b^i\text{ for each }i\geq 1, ab=a^2b^2\rangle$. Following the first suggestion, the author has shortened the proof that $S$ is combinatorial in Lemma \ref{22}(ii), originally proved by showing directly that $|H_{ba}|=1$. However, after examining the second suggestion, the author has concluded that it would not lead to a shorter way of describing combinatorial bisimple monogenic orthodox semigroups than that given in Section 2.
\med

\end{document}